\documentclass[sts]{imsart}

\usepackage{amssymb,amsbsy,amsmath,amscd,amsthm,amsfonts,MnSymbol}
\usepackage{cite}
\usepackage[utf8]{inputenc}
\usepackage{enumerate,hyperref}
\usepackage{dsfont}
\usepackage{graphicx}

\usepackage[top=1in,bottom=1in,right=1in,left=1in]{geometry}

\hypersetup{
hidelinks
}

\newtheorem{definition}{Definition}
\newtheorem{theorem}{Theorem}
\newtheorem{assumption}{Assumption}

\newtheorem{lemma}{Lemma}

\newtheorem{corollary}{Corollary}
\newtheorem{claim}{Proposition}

\newcommand{\R}{\mathbb R}
\newcommand{\PP}{\mathbb P}
\newcommand{\E}{\mathbb E}
\newcommand{\K}{\mathcal K_d}

\newcommand{\Sp}{\mathbb S}

\newcommand*\diff{\mathop{}\!\mathrm{d}}

\newcommand{\DS}{\displaystyle}

\date{}

\begin{document}

\begin{frontmatter}

\title{Uniform Behaviors of Random Polytopes under the Hausdorff Metric}
\runtitle{Random polytopes under the Hausdorff metric}

\begin{aug}

\author{\fnms{Victor-Emmanuel}~\snm{Brunel}\ead[label=veb]{veb@mit.edu}}

\affiliation{Massachusetts Institute of Technology}

\address{{Victor-Emmanuel Brunel}\\
{Department of Mathematics} \\
{Massachusetts Institute of Technology}\\
{77 Massachusetts Avenue,}\\
{Cambridge, MA 02139-4307, USA}\\
}

\runauthor{Brunel, V.-E.}
\end{aug}

\begin{abstract}

We study the Hausdorff distance between a random polytope, defined as the convex hull of i.i.d. random points, and the convex hull of the support of their distribution. As particular examples, we consider uniform distributions on convex bodies, densities that decay at a certain rate when approaching the boundary of a convex body, projections of uniform distributions on higher dimensional convex bodies and uniform distributions on the boundary of convex bodies. We essentially distinguish two types of convex bodies: those with a smooth boundary and polytopes. In the case of uniform distributions, we prove that, in some sense, the random polytope achieves its best statistical accuracy under the Hausdorff metric when the support has a smooth boundary and its worst statistical accuracy when the support is a polytope. This is somewhat surprising, since the exact opposite is true under the Nikodym metric. We prove rate optimality of most our results in a minimax sense. In the case of uniform distributions, we extend our results to a rescaled version of the Hausdorff metric. We also tackle the estimation of functionals of the support of a distribution such as its mean width and its diameter. Finally, we show that high dimensional random polytopes can be approximated with simple polyhedral representations that significantly decrease their computational complexity without affecting their statistical accuracy.

\end{abstract}


\end{frontmatter}

\section{Introduction}\label{Sec:Intro}

A simple representation of random polytopes consists of taking the convex hull of random points in the Euclidean space. Other representations have been suggested and studied, such as random projections of high dimensional polytopes or intersection of random halfspaces (about random polytopes, see \cite{Schneider2008, Reitzner2010, Hug2013} and the references therein).

The study of random polytopes goes back to R{\'e}nyi and Sulanke's seminal works \cite{RenyiSulanke1963, RenyiSulanke1964}, where the authors studied the area and the number of vertices of the convex hull of independent and identically distributed (i.i.d.) points uniformly distributed in a planar convex body. A long series of work followed so as to better understand the behavior of the volume of the random polytope in higher dimensions, depending on the structure of the supporting convex body \cite{Efron1965,SchneiderWieacker1980,BuchtaMuller1984,Dwyer1988,BaranyLarman1988,AffentrangerWieacker1991,Barany1992,BaranyBuchta1993,Schutt1994,Reitzner2003,Vu2005} also allowing for more general distributions (see\cite{Brunel2017} and the references therein).

Combinatorial properties of random polytopes are well understood in expectation \cite{Barany1989}
as well as very precise results about their functionals, such as intrinsic volumes
\cite{Muller1989,BoroczkyHoffmannHug2008,BoroczkyFodorReitznerVigh2009,BaranyFodorVigh2010}.

Beyond these probabilistic results, random polytopes are useful to approximate convex bodies by polytopes with few vertices \cite{Schneider1987,GlasauerSchneider1996,DumbgenWalther1996,Reitzner2004} and they are also in order to estimate density supports under convexity restrictions \cite{Moore1984,KTlectureNotes1993,KorostelevTsybakov1994,KorostelevSimarTsybakov1995,Cuevas2009,Brunel2016} (see also \cite{Pateiro2008} for weaker convexity type restrictions) or their volume \cite{Gayraud1997,BaldinReiss2016}.

In the aforementioned works, the random polytope is described as the convex hull of i.i.d. points whose distribution support is a convex body and an important quantity is its missing volume, i.e., the Lebesgue measure of the set-difference between the support and the random polytope. When the support has volume one, it is well known that the expected missing volume is of order at least $(\ln n)^{d-1}/n$ and at most $n^{-2/(d+1)}$, where $n$ is the number of points and $d$ is the dimension of the Euclidean space \cite{Groemer1974, Barany1989}. Moreover, the former rate is achieved when the support has a smooth boundary, whereas the latter is achieved when the support is a polytope \cite{Barany1989}. The Hausdorff distance between the random polytope and the convex support has attracted less interest and most significant results appear in works on approximation and estimation of convex bodies \cite{KorostelevTsybakov1994,DumbgenWalther1996}. In \cite{KorostelevTsybakov1994}, it is shown that the random polytope is a rate optimal estimator of the support, in a minimax sense, when the distribution of the random points is uniform. In \cite{DumbgenWalther1996}, more general distributions are considered, including uniform on the boundary of a convex body, or nearly uniform distributions in smooth convex bodies or in polytopes: Almost sure rates of convergence are proven for the Hausdorff distance between the support and the random polytope, as the number $n$ of points goes to infinity. One of the objectives of this article is to refine these results by proving general yet explicit deviation inequalities and show that the two extreme cases mentioned above (smooth boundary v.s. polytope) are reversed under the Hausdorff metric.

\section{Notation, outline and contributions}\label{Sec:Not}


Throughout this work, $d$ is a fixed positive integer that stands for the dimension of the ambient Euclidean space. A convex body in $\R^d$ is a compact and convex set with nonempty interior. We denote by $\mathcal K_d$ the class of all convex bodies in $\R^d$ and by $\mathcal K_d^{(1)}$ the class of all convex bodies included in the closed unit Euclidean ball in $\R^d$.

For all positive integers $p$, the $p$-dimensional closed Euclidean ball with center $a\in\R^p$ and radius $r>0$ is denoted by $B_p(a,r)$ and the $(p-1)$-dimensional Euclidean unit sphere is denoted by $\Sp^{p-1}$. We write $\kappa_p$ for the $p$-dimensional volume of $B_p(0,1)$ and $\omega_{p-1}$ for the surface area of the $(p-1)$-dimensional unit sphere (note that $\omega_{p-1}=p\kappa_p$). The $p$-dimensional volume of a Borel set $A$ in $\R^p$ is denoted by $\textsf{Vol}_p(A)$.

The Euclidean norm in $\R^d$ is denoted by $\|\cdot\|$ and the Euclidean distance is denoted by $\rho(\cdot,\cdot)$.

The interior of a set $A\subseteq \R^d$ is denoted by $\textsf{int}(A)$ and its boundary by $\partial A$. The complement of a set or event $A$ is denoted by $A^{\complement}$.

The support function $h_A$ of a compact set $A\subset \R^d$ is the mapping $\DS h_A(u)=\max\{\langle u,x\rangle:x\in A\}, u\in\R^d$, where $\langle\cdot,\cdot\rangle$ stands for the canonical inner product in $\R^d$: It is the largest signed distance between the origin and a supporting hyperplane of $A$ orthogonal to $u$. 

Support functions have two important properties, which will be useful in the sequel. Let $K$ be a compact and convex set. Then, $h_K$ is positively homogeneous, i.e.,
\begin{equation} \label{prop1}
	h_K(\lambda u)=\lambda h_K(u)
\end{equation}
for all $u\in\R^d, \lambda\geq 0$, and it is subadditive:
\begin{equation} \label{prop2}
	h_K(u+v) \leq h_K(u)+h_K(v)
\end{equation}
for all $u,v\in\R^d$. In particular, \eqref{prop2} implies the reverse triangle inequality
\begin{equation} \label{prop3}
	h_K(u-v) \geq h_K(u)-h_K(v)
\end{equation}
for all $u,v\in\R^d$. 

Let $K\in\K$ with $0\in\textsf{int}(K)$. The pseudo-norm $\|\cdot\|_K$ (also called the Minkowski functional of $K$) is defined as $\DS \|x\|_K=\min\{\lambda\geq 0 : x\in\lambda K\}, x\in\R^d$. It is a norm if and only if $K$ is symmetric, i.e., $K=-K$. If $K=B_d(0,1)$, then $\|\cdot\|_K$ is the Euclidean norm in $\R^d$. 

The polar set $K^\circ$ of a convex body $K$ is the set of all $x\in\R^d$ such that $h_K(x)\leq 1$. It is clear that $0$ is an interior point of $K^\circ$ as soon as $K$ is bounded. Moreover, if $0\in \textsf{int}(K)$, then $K^\circ$ is compact and $\DS h_K(x)=\|x\|_{K^\circ}$ for all $x\in\R^d$. In particular, if $K$ is bounded and has $0$ in its interior, then $K^{\circ}\in\K$.

The Hausdorff distance between two sets $K,K'\subseteq \R^d$ is denoted by $\textsf{d}_{\textsf{H}}(K,K')$ and it is defined as $\DS \textsf{d}_{\textsf{H}}(K,K')=\max\{\max_{x\in K}\rho(x,K'),\max_{y\in K'}\rho(y,K)\}$. If $K$ and $K'$ are compact and convex, it can be expressed in terms of their support functions:
$\DS \textsf{d}_{\textsf{H}}(K,K')=\sup_{u\in \Sp^{d-1}} |h_K(u)-h_{K'}(u)|$. 

If $K\in\mathcal K_d$, $u\in \Sp^{d-1}$ and $\varepsilon\geq 0$, we denote by $C_K(u,\varepsilon)$ the cap of $K$ in the direction of $u$ and of width $\varepsilon$:
$$C_K(u,\varepsilon)=\left\{x\in K: \langle u,x \rangle \geq h_K(u)-\varepsilon \right\}.$$

In the paper, all random variables are defined on a probability space $(\Omega,\mathcal F,\PP)$.  
For notation convenience and with no loss of generality, if $\mu$ is a given distribution in $\R^d$ equipped with its Borel $\sigma$-algebra, we identify $\PP$ and $\mu$ (by assuming that $\Omega=\R^d$). The support of $\mu$ is denoted by $\textsf{supp}(\mu)$. The uniform distribution on a compact set $A\subseteq \R^d$ is denoted by $\PP_A$ and the corresponding expectation operator is denoted by $\E_A$.

If $(u_n)$ and $(v_n)$ are two positive sequences, we write that $u_n=O(v_n)$ or, equivalently, $v_n=\omega(u_n)$, if the sequence $(u_n/v_n)$ is bounded.


The paper is organised as follows. In Section \ref{SectionMain}, we state and prove our main theorems. The first one is a uniform deviation inequality that yields a stochastic upper bound for the Hausdorff distance between the random polytope and the convex hull of the support of the random points. It is uniform in the sense that the bounds do not depend on the probability measure of the points or its support. The second one gives the rate of convergence of functionals of the random polytope, such as the mean width or the maximal width, to those of the convex hull of the support of the random points. Both theorems hold under an assumption on the probability mass of caps of the convex hull of the support of the random points. That assumption is discussed in the last part of Section \ref{SectionMain}, where we list particular cases that are most relevant to the stochastic geometry and the statistics literature. In Section \ref{Sec:Opt}, we assume that the distribution of the random points in uniform on a convex body and we prove rate-optimality of our first theorem in a minimax sense. In Section \ref{Sec:Pol}, we use a different technique in order to prove a refined deviation inequality for the convex hull of uniform points in polytopes. In addition, we show that up to a logarithmic factor, the polytopal case is the least favorable for the statistical accuracy of random polytopes, in some sense. This result is surprising since the opposite is true under the Nikodym metric (i.e., volume of the symmetric difference): polytopal supports are the most favorable, whereas convex supports with smooth boundary are the least \cite{Groemer1974, Barany1989}. In Section \ref{Sec:Resc}, we focus on uniform distributions and we extend our results under a rescaled version of the Hausdorff metric, which allows to consider a broader class of convex bodies. In Section \ref{Sec:CompComp}, we tackle the problem of computational complexity of high dimensional random polytopes and propose a random approximation of random polytopes that achieves a significant gain in the computational cost without affecting the statistical accuracy.

\section{Behavior of the random polytope and its functionals}\label{SectionMain}

\subsection{Behavior of the random polytope}

Let $K\in\K$ and $\mu$ be a probability measure in $\R^d$, whose support is contained in $K$. We make the following assumption on the pair $(\mu,K)$.

\begin{assumption}\label{Assump}
	Let $\alpha$, $L$ and $\varepsilon_0$ be positive numbers with $\varepsilon_0\leq 1$. For all $u\in\Sp^{d-1}$ and $\varepsilon\in [0,\varepsilon_0]$, $\DS \mu\left(C_K(u,\varepsilon)\right) \geq L\varepsilon^\alpha$. 
\end{assumption}

Note that if $(\mu,K)$ satisfy Assumption \ref{Assump}, then $K$ is necessarily the convex hull of the support of $\mu$. Indeed, $K$ is convex and it is easy to see that every closed halfspace that contains either $\textsf{supp}(\mu)$ or $K$ needs to contain the other as well. 

Consider a collection $X_1,\ldots,X_n$ of i.i.d. random points with probability distribution $\mu$, where $n$ is a positive integer and let $\hat K_n$ be their convex hull. The following theorem shows that under Assumption \ref{Assump}, $\hat K_n$ concentrates around $K$ under the Hausdorff metric, at an explicit rate that depends on $\alpha$. For all $\alpha>0$, we set 
$$C_\alpha=\inf_{t>0}\frac{(1+t)^\alpha}{1+t^\alpha}=\begin{cases} 1 \mbox{ if } \alpha\geq 1 \\ 2^{\alpha-1} \mbox{ otherwise.} \end{cases}$$

\begin{theorem}\label{theorem21}
Let $\alpha,L,\varepsilon_0$ be positive numbers with $0<\varepsilon_0\leq 1$. Assume Assumption \ref{Assump} and that $K\subseteq B_d(0,1)$. Set $\displaystyle{\tau_1=\max\left(1,\frac{d}{C_\alpha\alpha L}\right)}$, $\displaystyle{a_n=\left(\frac{\tau_1\ln n}{n}\right)^\frac{1}{\alpha}}$ and $\displaystyle{b_n=n^{\frac{-1}{\alpha}}}$. Then, the random polytope $\hat K_n$ satisfies 
	$$\PP\left[\textsf{d}_{\textsf{H}}(\hat K_n,K)\geq 2a_n + 2b_n x\right] \leq 12^d\exp\left(-C_\alpha L x^\alpha\right),$$
for all $x\geq 0$ with $\displaystyle{a_n+b_n x\leq \varepsilon_0}$.
\end{theorem}

\begin{proof}

Let $\varepsilon\in [0,\varepsilon_0]$ and $\delta=\varepsilon/4$.
For all $u\in\Sp^{d-1}$, we have $h_{\hat K_n}(u) \leq h_K(u)-\varepsilon$ if and only if all of the points $X_1,\ldots,X_n$ lie outside of $C_K(u,\varepsilon)$. Hence, 
\begin{align} \label{step1}
	\PP\left[h_{\hat K_n}(u)\leq h_K(u)-\varepsilon\right] & = \left(1-\mu\left(C_K(u,\varepsilon)\right)\right)^n \nonumber \\
	& \leq \exp\left(-Ln\varepsilon^\alpha\right).
\end{align}

The following lemma can be found in \cite[Lemma 5.2]{FresenVitale2014}. 

\begin{lemma} \label{lemma1}
Let $K\in\mathcal K_d$ with $0\in\textsf{int}(K)$ and $\delta\in (0,1/2)$. There exists a subset $\mathcal N_\delta\subseteq \partial K$ of cardinality at most $\displaystyle{\left(\frac{3}{\delta}\right)^d}$ such that, for all $u\in\partial K$, there exist two sequences $(u_j)_{j\geq 0}\subseteq \mathcal N_\delta$ and $(\delta_j)_{j\geq 1}\subseteq \R$ such that $0\leq \delta_j\leq \delta^j$ for all $j\geq 1$, and $\displaystyle{u=u_0+\sum_{j=1}^\infty \delta_j u_j}$.
\end{lemma}

Let $\mathcal N_\delta$ be a subset of $\Sp^{d-1}$ satisfying Lemma \ref{lemma1} applied to the unit Euclidean ball. Let us denote by $A$ the event $\{h_{\hat K_n}(u)> h_K(u)-\varepsilon, \forall u\in\mathcal N_\delta\}$. By \eqref{step1} and the union bound, 
\begin{equation} \label{step2}
	\PP(A)\geq 1- \left(\frac{3}{\delta}\right)^d\exp\left(-Ln\varepsilon^\alpha\right).
\end{equation}
Let $A$ hold and $u\in \Sp^{d-1}$. By Lemma \ref{lemma1}, we can write $\displaystyle{u=u_0+\sum_{j=1}^\infty \delta_j u_j}$ where $\displaystyle{(u_j)_{j\geq 0}\subseteq \mathcal N_\delta}$ and $0\leq\delta_j\leq \delta^j$ for all $j\geq 1$.
Note that almost surely, $\hat K_n\subseteq K\subseteq B_d(0,1)$, hence $h_{\hat K_n}(u)\leq h_K(u)\leq 1$, for all $u\in \Sp^{d-1}$. Thus,
\begin{align} \label{step3}
	h_{\hat K_n}(u) & = h_{\hat K_n}\left(u_0+\sum_{j=1}^\infty \delta_j u_j\right) \geq  h_{\hat K_n}(u_0)-\sum_{j=1}^\infty \delta^j h_{\hat K_n}(-u_j)\quad \mbox{ by \eqref{prop1} and \eqref{prop2}} \nonumber \\
	& > h_K(u_0)-\varepsilon-\frac{\delta}{1-\delta} \nonumber \\
	& \geq h_K(u)-\sum_{j=1}^\infty \delta^j h_K(u_j) - \varepsilon-\frac{\delta}{1-\delta} \quad \mbox{ again by \eqref{prop1} and \eqref{prop2}} \nonumber \\
	& \geq h_K(u)-\varepsilon-\frac{2\delta}{1-\delta} \geq h_K(u)-2\varepsilon,
\end{align}
since $\delta=\varepsilon/4 \leq \varepsilon_0/4\leq 1/2$.

Hence, by \eqref{step2} and \eqref{step3},
\begin{align} \label{step4}
	\PP\left[\textsf{d}_{\textsf{H}}(\hat K_n,K)\geq 2\varepsilon\right] & = \PP\left[\exists u\in \Sp^{d-1}, h_{\hat K_n}(u)\leq h_K(u)-\varepsilon\right] \nonumber \\
	& \leq \left(\frac{12}{\varepsilon}\right)^d\exp\left(-Ln\varepsilon^\alpha\right) \nonumber \\
	& \leq 12^d\exp\left(-Ln\varepsilon^\alpha -d\ln\varepsilon\right).
\end{align}
By setting $\varepsilon=a_n+b_n x$ for some nonnegative number $x$, where $a_n$ and $b_n$ are defined in Theorem \ref{theorem21}, we have that 
$$\varepsilon^\alpha\geq C_\alpha \left(\frac{\tau_1\ln n}{n} + \frac{x^\alpha}{n}\right),$$
and $$\ln \varepsilon \geq \ln a_n \geq -\frac{\ln n}{\alpha}.$$
The conclusion of Theorem \ref{theorem21} follows.

\end{proof}

If $\displaystyle{a_n+b_n x> 1}$, the upper bound in Theorem \ref{theorem21} can be replaced with $0$, since $K$ and $\hat K_n$ are both subsets of the unit ball, and so the Hausdorff distance between them is at most 2. When $\displaystyle{\varepsilon_0\leq a_n+b_n x\leq 1}$, the upper bound in Theorem \ref{theorem21} can be replaced with its corresponding value at $\displaystyle{x=(\varepsilon_0-a_n)b_n^{-1}}$, since the left-hand side of the inequality is a nonincreasing function of $x$. The constants in Theorem \ref{theorem21} do not depend on $\mu$ and $K$: The deviation inequality is uniform in all pairs $(\mu,K)$ that satisfy Assumption \ref{Assump}, with $K\subseteq B_d(0,1)$. It is important to notice that the constant factors in Theorem \ref{theorem21} are much smaller than in \cite[Theorem 1]{Brunel2017} when the dimension $d$ becomes large. However, we believe that the factor $12^d$ cannot be replaced with a sub-exponential factor in general.

Theorem \ref{theorem21} yields the following moment inequalities.

\begin{corollary}\label{corollary1}
	Let the assumptions of Theorem \ref{theorem21} hold. For all real numbers $q\geq 1$,
	$$\E\left[\textsf{d}_{\textsf{H}}(\hat K_n,K)^q\right]=O\left(\left(\frac{\ln n}{n}\right)^\frac{q}{\alpha}\right),$$
	with constant factors that do not depend on $\mu$ and $K$. 
\end{corollary}

\begin{proof}

Let $Z$ be a nonnegative random variable and $q$ be a positive number. Then, by Fubini's theorem,
$$\E[Z^q]=q\int_0^\infty t^{q-1}\PP[Z\geq t]dt.$$
Let us apply this equality to $Z=\textsf{d}_{\textsf{H}}(\hat K_n,K)$ (which we denote by $\textsf{d}_{\textsf{H}}$ for the sake of simplicity) and $q\geq 1$. Since $\textsf{d}_{\textsf{H}}(\hat K_n,K)\leq 2$ $\mu$-almost surely,
\begin{align}\label{step3331}
	\E[\textsf{d}_{\textsf{H}}&(\hat K_n,K)^q] & \nonumber \\
	& = q\int_0^\infty t^{q-1}\PP\left[\textsf{d}_{\textsf{H}}\geq t\right] dt  \nonumber \\
	& = q\int_0^{2a_n}  t^{q-1}\PP\left[\textsf{d}_{\textsf{H}}\geq t\right] dt + q\int_{2a_n}^{2\varepsilon_0} t^{q-1}\PP\left[\textsf{d}_{\textsf{H}}\geq t\right]dt \nonumber \\
	& \hspace{2cm} + q\int_{2\varepsilon_0}^2 t^{q-1}\PP\left[\textsf{d}_{\textsf{H}}\geq t\right]dt.
\end{align}
In the first integral, let us bound the probability by 1. In the second and third integrals, we perform the change of variable $t=2a_n +2b_nx$. Then, by Theorem \ref{theorem21} and the remark that follows it, \eqref{step3331} becomes:
\begin{align}\label{step3332}
	\E[\textsf{d}_{\textsf{H}}& (\hat K_n,K)^q] \nonumber \\
	& \leq (2a_n)^q + 2q b_n\int_0^{(\varepsilon_0-a_n)b_n^{-1}} (2a_n +2b_nx)^{q-1}\PP\left[\textsf{d}_{\textsf{H}}\geq 2a_n+2b_n x\right]dx \nonumber \\
	& \hspace{15mm} + 2q b_n\int_{(\varepsilon_0-a_n)b_n^{-1}}^{(1-a_n)b_n^{-1}} (2a_n +2b_nx)^{q-1}\PP\left[\textsf{d}_{\textsf{H}}\geq 2a_n+2b_n x\right]dx  \nonumber \\
	& \leq (2a_n)^q +  q b_n\int_0^{(\varepsilon_0-a_n)b_n^{-1}}(2a_n+2b_n x)^{q-1}12^d\exp\left(-C_\alpha L x^\alpha\right)dx \nonumber \\
	& \hspace{15mm} + q 2^{q-1}|1-\varepsilon_0|b_n^{-1}12^d\exp\left(-C_\alpha L |1-\varepsilon_0|^\alpha b_n^{-\alpha}\right). \nonumber \\
\end{align}
In the second term, we bound $(a_n +b_n x)^{q-1}$ by $2^{q-1}a_n^{q-1} +b_n^{q-1} x^{q-1}$, which yields
\begin{align} \label{step3333}
	& \int_0^{(\varepsilon_0-a_n)b_n^{-1}}(2a_n+2b_n x)^{q-1}12^d\exp\left(-C_\alpha L x^\alpha\right)dx \nonumber \\
	& \hspace{4mm}\leq 12^d 4^{q-1}\left(a_n^{q-1}\int_0^{\infty}\exp\left(-C_\alpha L x^\alpha\right)dx + b_n^{q-1}\int_0^{\infty}x^{q-1}\exp\left(-C_\alpha L x^\alpha\right)dx \right).
\end{align}
Since $b_n=O(a_n)$, \eqref{step3332} and \eqref{step3333} yield Corollary \ref{corollary1}. 

\end{proof}

Corollary \ref{corollary1} implies an upper bound for the minimax risk for estimation of $K\in\K^{(1)}$ under Assumption \ref{Assump}, of order $\DS \left((\ln n)/n\right)^{1/\alpha}$ (see \cite[Chapter 2]{Tsybakov2009} for a formal definition of minimax risks).

\subsection{Estimation of functionals of $K$}

In this section, we study a class of functionals of $K$ that can be written in terms of integrals involving its support function. We bound the performances of the corresponding functionals of $\hat K_n$, which are also called \textit{plug-in} estimators. What we call a functional is a mapping $T:\mathcal K_d\longrightarrow \R$ and the plug-in estimator of $T(K)$ is $T(\hat K_n)$.

For all real numbers $p\geq 1$ and all measurable functions $f:\Sp^{d-1}\rightarrow \R$, denote by 
$\DS \|f\|_p=\left(\int_{\Sp^{d-1}} |f(u)|^p d\sigma(u)\right)^\frac{1}{p}$,
where $\sigma$ is the uniform probability measure on the sphere and by $\|f\|_\infty=\sup_{u\in \Sp^{d-1}} |f(u)|$, provided these quantities are finite. 

For $K\in\mathcal K_d$, denote by $w_K(u)$ its width in direction $u\in\Sp^{d-1}$: $\DS w_K(u)= h_K(u)+h_K(-u)$. We are interested in the functionals $T_p(K)=\|h_K\|_p$ and $S_p(K)=\|w_K\|_p$, for $p\in[1,\infty]$ (extensions to broader classes of functionals would be possible). For instance, $S_1$ is the mean width and $S_\infty$ is the diameter. Asymptotics of the expected mean width of $\hat K_n$ are known in the case of uniform measures in $K$. When $\mu$ is uniform in a smooth convex body $K$, it is well known that $\DS S_1(K)-\E_K\left[S_1(\hat K_n)\right]\longrightarrow c_Kn^{\frac{-2}{d+1}}$ as $n\to\infty$, where $c_K$ is a positive number that depends on $K$, see \cite{SchneiderWieacker1980, BoroczkyFodorReitznerVigh2009} for explicit formulas. If $\mu$ is uniform in a polytope $K$, Schneider \cite{Schneider1987} proved that $\DS S_1(K)-\E_K\left[S_1(\hat K_n)\right] \longrightarrow c'_Kn^{\frac{-1}{d}}$ as $n\to\infty$,
for some other positive constant $c'_K$ that also depends on $K$. When $\mu$ is supported on the boundary of a smooth convex body $K$ with a positive density $g$ with respect to the surface area measure of $\partial K$, Müller \cite{Muller1989} showed that $\DS S_1(K)-\E_K\left[S_1(\hat K_n)\right] \longrightarrow c''_{g,K}n^{\frac{-2}{d-1}}$ as $n\to\infty$,
where the positive number $c''_{g,K}$ depends on $g$ and $K$. Here, we prove moment inequalities for $S_1(\hat K_n)$ as well as for the other functionals $S_p, T_p, p\in[1,\infty]$ under Assumption \ref{Assump}. 

Let $p\in [1,\infty]$ and $q\geq 1$. If $K$ and $K'$ are two convex bodies, then by the triangle inequality,
\begin{equation} \label{TriangleIneqT}
	\left|T_p(K)-T_p(K')\right|^q \leq \|h_K-h_{K'}\|_p^q
\end{equation}
and
\begin{equation} \label{TriangleIneqS}
	\left|S_p(K)-S_p(K')\right|^q \leq \|w_K-w_{K'}\|_p^q.
\end{equation}
Next theorem states uniform upper bounds for the moments $\displaystyle{\E\left[\|h_K-h_{\hat K_n}\|_p^q\right]}$ and $\displaystyle{\E\left[\|w_K-w_{\hat K_n}\|_p^q\right]}$ when $(\mu,K)$ satisfies Assumption \ref{Assump}. Note that the case $p=\infty$ is treated in Corollary \ref{corollary1} and we only treat here the case of finite $p$.
\begin{theorem}\label{theorem2'}
	Let $d$ be a positive integer, $\alpha,L>0$, $0<\varepsilon_0\leq 1$ and $p,q\geq 1$ be real numbers. Assume Assumption \ref{Assump} and that $K\subseteq B_d(0,1)$. Then,
	$$\E\left[\|h_K-h_{\hat K_n}\|_p^q\right]=O\left(n^{\frac{-q}{\alpha}}\right)$$
	and 
	$$\E\left[\|w_K-w_{\hat K_n}\|_p^q\right]=O\left(n^{\frac{-q}{\alpha}}\right).$$
\end{theorem}

\begin{proof}

Let $u\in \Sp^{d-1}$. By \eqref{step1}, for all $\varepsilon\in [0,\varepsilon_0]$,
$$\PP\left[h_{\hat K_n}(u)\leq h_K(u)-\varepsilon\right] \leq \exp\left(-Ln\varepsilon^\alpha\right).$$
For $\varepsilon\in (\varepsilon_0,2]$,
$$\PP\left[h_{\hat K_n}(u)\leq h_K(u)-\varepsilon\right] \leq \exp\left(-Ln\varepsilon_0^\alpha \right),$$
since the left hand side is a nonincreasing function of $\varepsilon$ and if $\varepsilon > 2$,
$$\PP\left[h_{\hat K_n}(u)\leq h_K(u)-\varepsilon\right]=0.$$
Hence, by a similar argument as in the proof of Corollary \ref{corollary1}, for all real numbers $k\geq 1$, there is some positive number $c_k$ that depends neither on $\mu, K$ nor on $u$ and such that
\begin{equation}\label{step2'1}
	\E\left[|h_K(u)-h_{\hat K_n}(u)|^k\right] \leq c_kn^{\frac{-k}{\alpha}}.
\end{equation}
If $q\leq p$, Jensen's inequality applied to the concave function $\displaystyle{[0,\infty)\ni x\mapsto x^\frac{q}{p}}$ yields
\begin{align} \label{step2'1bis}
	\E\left[\|h_{\hat K_n}-h_K\|_p^q\right] & = \E\left[\left(\int_{\Sp^{d-1}}|h_K(u)-h_{\hat K_n}(u)|^pd\sigma(u)\right)^{q/p}\right] \nonumber \\
	& \leq \E\left[\int_{\Sp^{d-1}}|h_K(u)-h_{\hat K_n}(u)|^pd\sigma(u)\right]^{q/p}.
\end{align}
By Fubini's theorem, which allows switching the expectation and the integral and by \eqref{step2'1} with $k=p$, \eqref{step2'1bis} yields
\begin{equation*}
	\E\left[\|h_{\hat K_n}-h_K\|_p^q\right] \leq c_p^{q/p}n^{\frac{-q}{\alpha}},
\end{equation*}
which proves the theorem in the case $q\leq p$. If now $q>p$, the mapping $\displaystyle{[0,\infty)\ni x\mapsto x^\frac{q}{p}}$ is convex and Jensen's inequality yields
\begin{equation*}
	\left(\int_{\Sp^{d-1}}|h_K(u)-h_{\hat K_n}(u)|^p\right)^{q/p} \leq \int_{\Sp^{d-1}}|h_K(u)-h_{\hat K_n}(u)|^q,
\end{equation*}
so again by Fubini's theorem and using \eqref{step2'1} with $k=q$,
\begin{equation*}
	\E\left[\|h_{\hat K_n}-h_K\|_p^q\right] \leq c_qn^{\frac{-q}{\alpha}}.
\end{equation*}
\end{proof}

It is not surprising that in Theorem \ref{theorem2'}, there is no logarithmic factor when $p\neq\infty$, as opposite to the case $p=\infty$ (see Corollary \ref{corollary1}). This is due to the fact that $L^\infty$ norms are more restrictive than $L^p$ norms ($p\in[1,\infty)$), since they capture the largest values of a function, no matter the size of the subset where large values are achieved. When $\mu$ is the uniform probability measure in a smooth convex body, Theorem \ref{theorem3} in Section \ref{Sec:Opt} shows that the logarithmic factor cannot be avoided for $p=\infty$. In \cite{Guntuboyina2012}, a convex body is estimated from noisy observations of its support function, and the risk is evaluated in terms of the $L^2$ norm of the support function, which exactly corresponds to $p=2$ in Theorem \ref{theorem2'}. The risk of the least square estimator does not have logarithmic factors, and it is proven to be rate optimal, in a minimax sense. However, we believe that a logarithmic factor would be unavoidable if the risk was measured in terms of the $L^\infty$ norm of the support function, i.e., in terms of the Hausdorff metric. 
%
%
%

Another functional of interest, the \textit{thickness} of $K$, exhibits a similar behavior as the diameter. The thickness of $K$, denoted by $\Delta(K)$, is its minimal width, i.e., $\DS \Delta(K)=\min_{u\in\Sp^{d-1}} w_K(u)$. It is easy to see that $\DS |\Delta(\hat K_n)-\Delta(K)|\leq \|w_{\hat K_n}-w_K\|_{\infty}$. Hence, the plug-in estimator $\Delta(\hat K_n)$ converges at least at the same speed as $\hat K_n$ converges to $K$ under the Hausdorff metric.

\subsection{Special cases of Assumption \ref{Assump}} \label{Subsec:Spec}

Before listing special cases, we recall the definition of the \textit{reach} of a set (see, e.g., \cite[Definition 11]{Thale2008}).

\begin{definition}
	Let $A\subset\R^d$. The \textit{reach} of $A$ is the supremum of the set of real numbers $r\geq 0$ such that each $x\in A^{\complement}$ with $\rho(x,\partial A)\leq r$ has a unique metric projection on $\partial A$.
\end{definition}

A geometric interpretation of the reach is the radius of the smallest Euclidean ball that can be rolled up on the boundary of $A$, on the outside of $A$. If $A$ has reach at least $r>0$, we also say that the complement of $A$ satisfies the \textit{$r$-rolling ball condition}. In particular, if $K$ is a convex body and its complement has reach $r>0$, then for each point $x\in\partial K$, there exist $a\in K$ such that $x\in B_d(a,r)\subseteq K$. 

For all real numbers $r\in (0,1]$, we denote by $\mathcal K_{d,r}$ the class of all convex bodies in $\R^d$ whose complement has reach at least $r$ and by $\mathcal K_{d,r}^{(1)}$ the class of all convex bodies included in $B_d(0,1)$ whose complement has reach at least $r$.

\paragraph{Nearly uniform distributions in a convex body with positive reach complement}

Let $K\in\mathcal K_{d,r}^{(1)}$ and $\mu$ be a probability measure with $\textsf{supp}(\mu)=K$ and $\mu(A)\geq \lambda \textsf{Vol}_d(A)$, for all Borel sets $A\subseteq K$, where $\lambda$ is some given positive number. If $\mu$ is the uniform probability measure on $K$, take $\lambda=\kappa_d^{-1}$. 

\begin{claim} \label{Claim1}
The pair $(\mu,K)$ satisfies Assumption \ref{Assump} with $\alpha=(d+1)/2$, $\DS L=\frac{2\lambda\kappa_{d-1}r^\frac{d-1}{2}}{d+1}$ and $\varepsilon_0=r$.
\end{claim}

\begin{proof}
Let $u\in\Sp^{d-1}$ and $x^*\in \partial K$ such that $h_K(u)=\langle u,x^*\rangle$. Since $K$ satisfies the $r$-rolling ball property, $x^*\in B_d(a,r)\subseteq K$ for some $a\in K$. Let $B=B_d(a,r)$. Then, for all $\varepsilon\in [0,r]$,
\begin{align*}
	\mu\left(C_K(u,\varepsilon)\right) & \geq \lambda \textsf{Vol}_d\left(C_K(u,\varepsilon)\right) \\
	& \geq \lambda \textsf{Vol}_d\left(C_B(u,\varepsilon)\right) \\
	& = \lambda\int_0^\varepsilon \left(x(2r-x)\right)^\frac{d-1}{2}\kappa_{d-1} dx \\
	& \geq \lambda\kappa_{d-1}r^\frac{d-1}{2}\int_0^\varepsilon x^\frac{d-1}{2} dx \\
	& = \frac{2\lambda\kappa_{d-1}r^\frac{d-1}{2}}{d+1}\varepsilon^\frac{d+1}{2}.
\end{align*}
\end{proof}

\paragraph{Projection of a high dimensional uniform distribution in a convex body}

Let $D>d$ be an integer and $\tilde K\in\mathcal K_{D,r}^{(1)}$. Denote by $\pi_d$ the orthogonal projection in $\R^D$ onto the first $d$ coordinates and let $K=\pi_d(\tilde K)$. Then, if we identify $\R^d$ to $\R^d\times\{0\}^{D-d}$, it holds that $K\in \mathcal K_{d,r}^{(1)}$. Let $\mu$ be the image of the uniform distribution on $\tilde K$ by $\pi_d$, i.e., the distribution of $\pi_d(\tilde X)$, where $\tilde X$ is a uniform random variable in $\tilde K$.

\begin{claim} \label{Claim:Proj}
The pair $(\mu,K)$ satisfies Assumption \ref{Assump} with $\alpha=(D+1)/2$, $\DS L=\frac{2\kappa_{D-1}r^{\frac{D-1}{2}}}{\kappa_D(D+1)}$ and $\varepsilon_0=r$.
\end{claim}

\begin{proof}
This proposition follows from a similar argument as in the previous case. Indeed, let $u\in\Sp^{d-1}$, $x^*\in\partial K$ with $h_K(u)=\langle u,x^*\rangle$ and $a\in K$ such that $x^*\in B_d(a,r)\subseteq K$. Let $\tilde u$ be the (unique) unit vector in $R^D$ such that $u=\pi_d(\tilde u)$. Let $\tilde x^*\in\partial \tilde K$ with $h_{\tilde K}(\tilde u)=\langle \tilde u,\tilde x^*\rangle$ and $\tilde a\in \tilde K$ such that $\tilde x^*\in B_D(\tilde a,r)\subseteq \tilde K$. Denote by $\tilde \mu$ the uniform distribution in $\tilde K$. Then, for all $\varepsilon\in [0,r]$, $\DS \mu\left(C_K(u,\varepsilon)\right) = \tilde \mu\left(C_{\tilde K}(\tilde u,\varepsilon)\right)$ and one can use the same argument as in the end of the proof of Proposition \ref{Claim1}.
\end{proof}

More generally, one can consider distributions with densities that have a polynomial decay near the boundary of $K$.

\paragraph{Densities with polynomial decay near the boundary of $K$}

Let $K\in\mathcal K_{d,r}^{(1)}$ and $\mu$ be a probability measure with $\textsf{supp}(\mu)=K$. Assume that $\mu$ has a density $f$ with respect to the Lebesgue measure and that for all $x\in K$ with $\rho(x,\partial K)\leq r$, $f(x)\geq C\rho(x,\partial K)^\gamma$, for some real numbers $C>0$ and $\gamma\geq 0$. For instance, this is satisfied by the previous case with $\DS\gamma=\frac{D-d}{2}$ and $\DS C=r^{\frac{D-d}{2}}\kappa_{D-d}$. Note that the exponent $\alpha$ in Proposition \ref{Claim:Proj}, which we proved using geometric arguments, could also be deduced from this fact together with the following proposition. 

\begin{claim}
	The pair $(\mu,K)$ satisfies Assumption \ref{Assump} where $\DS \alpha=\frac{d+1+2\gamma}{2}$, $L>0$ is a positive constant that only depends on $r$, $d$ and $\gamma$ and $\varepsilon_0=r$.
\end{claim}

The explicit form of the constant $L$ can be found in the following computation. 

\begin{proof}
Let $u\in\Sp^{d-1}$ and $\varepsilon\in [0,r)$. Let $x^*\in\partial K$ with $h_K(u)=\langle u,x^*\rangle$ and $B=B_d(a,r)$ a ball of radius $r$ containing $x^*$ and included in $K$. Then, 
\begin{align*}
	\mu(C_K(u,\varepsilon))& \geq \mu(C_B(u,\varepsilon)) \geq C\int_{C_B(u,\varepsilon)}\rho(x,\partial K)^\gamma \diff x \geq C\int_{C_B(u,\varepsilon)}\rho(x,\partial B)^\gamma \diff x \\
	& = C(d-2)\omega_{d-2}\int_0^\varepsilon\left(\int_0^{\sqrt{2rt-t^2}}s^{d-2}\left(r-\sqrt{s^2+(r-t)^2}\right)^\gamma\diff s\right)\diff t \\
	& \geq \frac{C(d-2)\omega_{d-2}}{2r}\int_0^\varepsilon\left(\int_0^As^{d-2}\left(A^2-s^2\right)^\gamma\diff s\right)\diff t \\
	& \geq \frac{C(d-2)\omega_{d-2}\mathcal B(d-1,\gamma+1)}{2r}\int_0^\varepsilon \left(\sqrt{2rt-t^2}\right)^{2\gamma+d-1}\diff t \\
	& \geq \frac{C(d-2)\omega_{d-2}\mathcal  B(d-1,\gamma+1)r^{(2\gamma+d-3)/2}}{2}\varepsilon^{(2\gamma+d+1)/2},
\end{align*}
where we identified the unit ball in the hyperplane that is orthogonal to $u$ with $B_{d-1}(0,1)$, we denoted by $A=\sqrt{2rt-t^2}$ and by $\mathcal B(\cdot,\cdot)$ the Beta function. 
\end{proof}

\paragraph{Nearly uniform distributions on the boundary of a convex body with positive reach complement}

Let $K\in\mathcal K_{d,r}^{(1)}$ and $\mu$ be a probability measure supported on $\partial K$ with density $f$ with respect to the $(d-1)$-dimensional Hausdorff measure on $\partial K$, satisfying $f(x)\geq\lambda$, for all $x\in \partial K$ and some $\lambda>0$. Note that if $\mu$ is the uniform distribution on $\partial K$, then one can take $\lambda=1/\omega_{d-1}$, since $K\subseteq B_d(0,1)$. 

\begin{claim}
	The pair $(\mu, K)$ satisfies Assumption \ref{Assump} with $\alpha=(d-1)/2$, $L=\lambda \omega_{d-1}r^{(d-1)/2}$ and $\varepsilon_0=r$.
\end{claim}

\begin{proof}
Let $u\in \Sp^{d-1}$ and $0\leq\varepsilon\leq r$. Then, $\mu(C_K(u,\varepsilon))\geq \lambda \mathcal H_{d-1}(\partial K\cap C_K(u,\varepsilon))$, where $\mathcal H_{d-1}$ stands for the $(d-1)$-dimensional Hausdorff measure. Since $K$ satisfies the $r$-rolling condition, $\mathcal H_{d-1}(\partial K\cap C_K(u,\varepsilon))\geq \mathcal H_{d-1}(\partial B\cap C_B(u,\varepsilon))$, where $B$ is any $d$-dimensional Euclidean ball with radius $r$. It follows that
\begin{align*}
	\mathcal H_{d-1}(\partial K\cap C_K(u,\varepsilon)) & \geq \mathcal H_{d-1}(\partial B\cap C_B(u,\varepsilon)) \\
	& \geq \frac{1}{2}\omega_{d-1}r^{d-1}\int_0^{\frac{\varepsilon(2r-\varepsilon)}{r^2}}t^\frac{d-3}{2}(1-t)^{-1/2} dt \\
	& \geq \frac{1}{2}\omega_{d-1}r^{d-1}\int_0^{\frac{\varepsilon}{r}}t^\frac{d-3}{2} dt \\
	& \geq \omega_{d-1}r^\frac{d-1}{2}\varepsilon^\frac{d-1}{2}
\end{align*}
In addition, since $K$ is convex and is included in $B_d(0,1)$, the surface area of its boundary is bounded from above by that of $\Sp^{d-1}$, i.e., $\omega_{d-1}$. 
\end{proof}

\paragraph{Nearly uniform distributions in a convex body with general boundary}

Let $K\in\mathcal K_d^{(1)}$ and let $a\in\R^d$ and $\eta>0$ such that $B_d(a,\eta)\subseteq K$. The existence of $a$ and $\eta$ is guaranteed by the fact that the convex body $K$ must have nonempty interior. Let $\mu$ be a probability measure with $\textsf{supp}(\mu)=K$ and that satisfies $\mu(A)\geq \lambda \textsf{Vol}_d(A)$ for all Borel sets $A\subseteq K$, where $\lambda>0$.

\begin{claim} \label{Claim:Gen}
	The pair $(\mu,K)$ satisfies Assumption \ref{Assump} with $\alpha=d$, $\DS L=\frac{\lambda\kappa_{d-1}}{3\eta^{d-1}}$ and $\varepsilon_0=\eta$.
\end{claim}

\begin{figure} \label{Fig:GenPol}
\centering
   \includegraphics[scale=0.30]{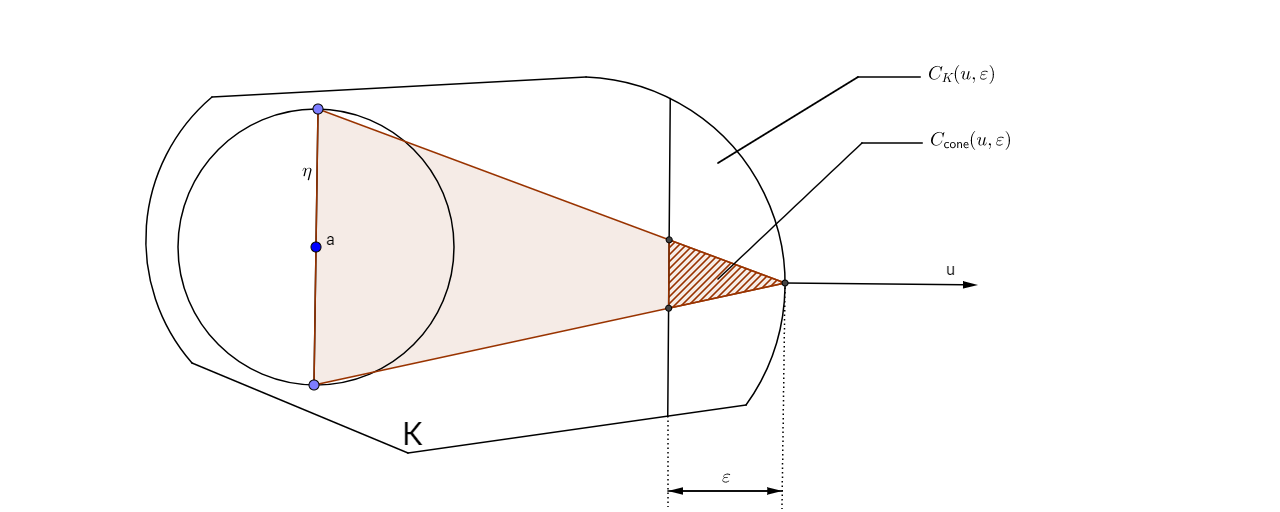}
   \caption{\label{etiquette2} The case of general convex bodies}
\end{figure}

\begin{proof}
Let $u\in\Sp^{d-1}$ and let $x^*\in\partial K$ such that $h_K(u)=\langle u,x^*\rangle$. Then, for all $\varepsilon\in [0,\eta]$, $C_K(u,\varepsilon)\supseteq C_{\textsf{cone}}(u,\varepsilon)$, where $\textsf{cone}$ is the smallest cone with apex $x^*$ containing the $(d-1)$-dimensional ball with center $a$, radius $\eta$, in the affine hyperplane containing $a$, orthogonal to $u$ (cf. Figure \ref{etiquette2}). Hence, $\DS \mu(C_K(u,\varepsilon))\geq \mu(C_{\textsf{cone}}(u,\varepsilon))\geq \frac{\lambda \varepsilon^d\kappa_{d-1}}{3\eta^{d-1}}$.

\end{proof}

\section{Rate-optimality of the random polytope}\label{Sec:Opt}

Let $n$ be a positive integer and $X_1,\ldots,X_n$ be random variables mapping the probability space $(\Omega,\mathcal F,\PP)$ onto $\R^D$, for some integer $D\geq 1$. A set-valued statistic is a set of the form $\tilde K_n(X_1,\ldots,X_n)$, where $\tilde K_n$ maps $\DS \left(\R^D\right)^n$ onto the collection of all compact subsets of $\R^d$, such that the function that maps $(\omega,x)\in \Omega\times \R^d$ to 1 if $\DS x\in \tilde K_n\left(X_1(\omega),\ldots,X_n(\omega)\right)$ and is zero otherwise is measurable with respect to the product of $\mathcal F$ and the Borel $\sigma$-algebra of $\R^d$. For the sake of brevity, we write $\tilde K_n$ instead of $\tilde K_n(X_1,\ldots,X_n)$. For instance, if $X_1,\ldots,X_n$ are i.i.d. random points, then $\hat K_n$ is a set-valued statistic. If a set-valued statistic $\tilde K_n$ is aimed to estimate a compact set $K$ which may be unknown in practice, it is called an estimator of $K$ and $\textsf{d}_{\textsf{H}}(\tilde K_n,K)$ is a (measurable) random variable. More general definitions of set-valued random variables can be found in \cite{Molchanov2005}, but they are not needed for our purposes.

\paragraph{Uniform distribution in a convex body with positive reach complement}

Next theorem, together with Corollary \ref{corollary1}, establishes rate minimax optimality of $\hat K_n$ on the class $\mathcal K_{d,r}^{(1)}$, under the Hausdorff metric.

\begin{theorem}\label{theorem3}
Let $r\in (0,1)$. For all real numbers $q\geq 1$, 
$$\inf_{\tilde K_n}\sup_{K\in\mathcal K_{d,r}^{(1)}} \E_K\left[\textsf{d}_{\textsf{H}}(\tilde K_n,K)^q\right] = \omega\left(\left(\frac{\ln n}{n}\right)^\frac{2q}{d+1}\right),$$
where the infimum is taken over all set-valued statistics $\tilde K_n$.
\end{theorem}

Note that for $r=1$, the class $\mathcal K_{d,1}^{(1)}$ reduces to the singleton $\left\{B_d(0,1)\right\}$ and the trivial estimator $\tilde K_n=B_d(0,1)$ satisfies $\displaystyle{\sup_{K\in\mathcal K_{1,d}^{(1)}} \E_K\left[\textsf{d}_{\textsf{H}}(\tilde K_n,K)\right]=0}$.

\begin{proof}

Here, we extend the arguments of the proof of \cite[Theorem 3]{KorostelevTsybakov1994} to higher dimensions. This standard method was also exploited in \cite{KorostelevSimarTsybakov1995'}.
We fix $R\in(r,1]$ and set $G_0=B_d(0,R)$. Let $\delta$ be a positive and small enough real number. For $u\in \Sp^{d-1}$, we define the set $G(u)$ as follows.
Let $\eta:\R\rightarrow\R$ be a given nonnegative function, twice continuously differentiable, such that $\eta(x)=0$ for all $x$ of norm at least one and $\DS \max_{x\in\R} \eta(x)=\eta(0)=1$. For instance, one can take $\displaystyle{\eta(x)=e^4g(2x-1)g(2-2x)}$, for all $x\in\R$, where $g(x)=\exp(-1/x)\mathds 1_{x>0}$. Let $H$ be the supporting hyperplane of $G$ at the point $Ru\in\partial G_0$. Identify this hyperplane with $\R^{d-1}$, with origin at the point $Ru$. Then, a parametrization of $G_0$ is 
$$G_0=\left\{(t,y)\in H\times\R :\|t\|\leq R, \hspace{2mm}R-\sqrt{R^2-\|t\|^2}\leq y\leq R+\sqrt{R^2-\|t\|^2}\right\}.$$
We define the set $G(u)$ by modifying the parametrisation of $\partial G_0$ in a small neighborhood of $Ru$:
\begin{align*}
	G(u)=\Big\{(t,y)\in H\times \R : & \|t\|\leq R, \\
	& R-\sqrt{R^2-\|t\|^2}+\alpha\delta^2\eta\left(\frac{2\|t\|}{R\delta}\right) \leq y\leq R+\sqrt{R^2-\|t\|^2}\Big\},
\end{align*}
where $\alpha$ is a positive number that can be tuned inependently of $\delta$ so that $G(u)\in\mathcal K_{d,r}$, for small enough $\delta$. Note that for all $u\in \Sp^{d-1}, G(u)\subseteq G_0$ and
\begin{align}\label{step31}
	\textsf{Vol}_d\left(G_0\backslash G(u)\right) & = \alpha\delta^2\int_{\R^d}\eta\left(\frac{2\|t\|}{R\delta}\right)dt \nonumber \\
	& = \frac{\alpha\delta^{d+1}R^{d-1}}{2^{d-1}}\int_{\R^d}\eta(\|t\|)dt \nonumber \\
	& = c\delta^{d+1},
\end{align}
for some positive constant $c$ that depends on $R$ and $d$ only. In addition, if $(u,u')\in \Sp^{d-1}\times \Sp^{d-1}$ is such that $\|u-u'\|\geq\delta$, 
\begin{equation} \label{step32}
	\textsf{d}_{\textsf{H}}\left(G(u),G(u')\right)=\alpha\delta^2.
\end{equation}
Let now $\mathcal N$ be a maximal $\delta$-packing of $\Sp^{d-1}$, i.e., a finite subset of $\Sp^{d-1}$ with maximal cardinality, satisfying that for all $(u,v)\in\mathcal N$ with $u\neq v$, $\|u-v\| \geq\delta$. Let us denote by $N$ its cardinality and by $u_1,\ldots,u_N$ its elements and, for $j=1,\ldots,N$, set $G_j=G(u_j)$. By a standard argument (e.g. \cite{KolomogorovTikhomirov1959}), $N\geq c\delta^{-d}$, for some positive constant $c$.
Now, by setting $\displaystyle{\delta=\left(\frac{\ln n}{n}\right)^\frac{1}{d+1}}$, we apply \cite[Lemma 1]{KorostelevSimarTsybakov1995'} to the sets $G_0,G_1,\ldots,G_N$, whose assumptions are satisfied due to \eqref{step31} and \eqref{step32} and this proves Theorem \ref{theorem3} for $q=1$. The lower bound for larger values of $q$ is then obtained by H\"older's inequality.

\end{proof}

\paragraph{Projection of a high dimensional uniform distribution in a convex body with positive reach complement}

The same idea as for the previous proof can be carried over in order to prove that $n^{-2/(D+1)}$ is the minimax rate of estimation of the projection of a $D$-dimensional convex support with positive reach ($D>d$).

Let $D>d$ be an integer. Let $\tilde K\in\mathcal K_D$ and $K=\pi_d(\tilde K)$ be the orthogonal projection of $\tilde K$ on the first $d$ coordinates. Let $\tilde X_1,\ldots,\tilde X_n$ be i.i.d. uniform random variables in $\tilde K$ and $X_1,\ldots,X_n$ be their orthogonal projections on the first $d$ coordinates. We denote by $\E_{\tilde K}^{\textsf{proj}_d}$ the expectation operator associated with the distribution of $X_1,\ldots,X_n$.

\begin{theorem}\label{Thm:OptProj}

Let $D>d$ be a positive integer and $0<r<1$. Then, for all real numbers $q\geq 1$,
$$\inf_{\tilde K_n}\sup_{\tilde K\in\mathcal K_{d,r}^{(1)}} \E_{\tilde K}^{\textsf{proj}_d}\left[\textsf{d}_{\textsf{H}}(\tilde K_n,\pi_d(\tilde K))^q\right] = \omega\left(\left(\frac{\ln n}{n}\right)^\frac{2q}{D+1}\right),$$
where the infimum is taken over all set-valued statistics $\tilde K_n$. 

\end{theorem}

\begin{proof}

Let $R\in (r,1]$ and $\delta>0$. Consider a $\delta$-packing $\mathcal N$ of the $(d-1)$-dimensional sphere $\Sp^{d-1}\times \{0\}^{D-d}$. Let $\tilde K=B_D(0,R)$ and for $u\in\mathcal N$, let $\tilde K(u)$ be defined exactly as $G(u)$ in the proof of Theorem \ref{theorem3}. Let $K=\pi_d(\tilde K)$ and $K(u)=\pi_d(\tilde K(u))$ for all $u\in\mathcal N$. Then, similarly to \eqref{step31} and \eqref{step32}, $\textsf{Vol}_d\left(\tilde K\setminus \tilde K(u)\right)=c\delta^{D+1}$ and $\textsf{d}_{\textsf{H}}(K(u),K(u'))=\alpha\delta^2$ for all $u,u'\in\mathcal N$ with $u\neq u'$. Denote by $\textsf{KL}$ the Kullback-Leibler divergence between probability measures (cf. \cite[Definition 2.5]{Tsybakov2009}). Then, by the data processing inequality (cf. \cite[Lemma 2]{Stoltz2017}), for all $u\in \mathcal N$,
\begin{align*}
\textsf{KL}(\pi_d\PP_{\tilde K},\pi_d\PP_{\tilde K(u)}) & \leq \textsf{KL}(\PP_K,\PP_{\tilde K(u)}) = \ln\left(\frac{\textsf{Vol}_d\left(\tilde K\right)}{\textsf{Vol}_d\left(\tilde K(u)\right)}\right) \\
& \leq \frac{\textsf{Vol}_d\left(\tilde K\setminus \tilde K(u)\right)}{\textsf{Vol}_d\left(\tilde K(u)\right)}\leq C\delta^{D+1},
\end{align*}
for some positive constant $C$ that depends on $d$ only and where we denoted by $\pi_d\PP_{\tilde K}$ (resp. $\pi_d\PP_{\tilde K(u)}$) the image of the distribution $\PP_{\tilde K}$ (resp. $\PP_{\tilde K(u)}$) by the orthogonal projection $\pi_d$. Then, Fano's lemma (cf. \cite[Corollary 2.6]{Tsybakov2009}) yields Theorem \ref{Thm:OptProj} for $q=1$, and the result can be obtained for larger values of $q$ by using H\"older's inequality.

\end{proof}

\paragraph{Uniform distribution in a general convex body}

Here, we prove that in the case of a uniform distribution in some convex body $K\in\K$ and under no assumption on the boundary of $K$, the rate that we achieved for polytopes in Section \ref{Sec:Pol}, i.e., $n^{-1/d}$, cannot be improved in a minimax sense. 

\begin{theorem}\label{Thm:OptPol}

For all real numbers $q\geq 1$,
$$\inf_{\tilde K_n}\sup_{K\in\mathcal K_d^{(1)}} \E_K\left[\textsf{d}_{\textsf{H}}(\tilde K_n,K)^q\right] = \omega\left(n^{-q/d}\right),$$
where the infimum is taken over all set-valued statistics $\tilde K_n$.
\end{theorem}

\begin{proof}
This proof is an extension of the proof of \cite[Theorem 2]{KorostelevTsybakov1994}. Let $K_0=\{x=(x_1,\ldots,x_d)\in B_d(0,1):x_j\geq 0,\forall j=1,\ldots,d\}$ and $K_1=\{x=(x_1,\ldots,x_d)\in B_d(0,1):x_j\geq 0,\forall j=1,\ldots,d, x_1+\ldots+x_d\geq h\}$, where $h>0$. It is easy to see that $\textsf{d}_{\textsf{H}}(K_0,K_1)=h$ and $\textsf{Vol}_d\left(G_0\triangle G_1\right)=ch^d$, where $c$ is the volume of the $d$-dimensional simplex. The same argument as in \cite[Theorem 2]{KorostelevTsybakov1994} yields Theorem \ref{Thm:OptPol}. 
\end{proof}

\section{Uniform distributions with polytopal supports} \label{Sec:Pol}

Let $K$ be a polytope and denote by $V$ its vertex set. Next lemma follows from the fact that a convex function defined on a compact and convex set achieves its maximum at an extreme point of that convex set.

\begin{lemma}
	Let $L$ be a convex set included in $K$. Then, 
	$$\textsf{d}_{\textsf{H}}(K,L)=\max_{v\in V}\rho(v,L).$$
\end{lemma}

It follows that, for $\varepsilon>0$, $\textsf{d}_{\textsf{H}}(K,\hat K_n)\geq\varepsilon$ if and only if $X_i\notin B_d(v,\varepsilon)$, for some vertex $v\in V$, and for all $i=1,\ldots,n$. Hence, by a union bound, if $\textsf{Vol}_d(K)\neq 0$,
\begin{align} \label{PolytoCase}
	\PP_K\left[\textsf{d}_{\textsf{H}}(K,\hat K_n)\geq \varepsilon\right] & \leq \sum_{v\in V} \PP_K\left[X_i\notin B_d(v,\varepsilon), \forall i=1,\ldots,n\right] \nonumber \\
	& \leq \sum_{v\in V} \left(1-\frac{\textsf{Vol}_d\left(B_d(v,\varepsilon)\cap K\right)}{\textsf{Vol}_d(K)}\right)^n.
\end{align}
 
In order to get a uniform inequality on a certain class of polytopes, it is natural to restrict the polytope $K$ to have a bounded number a vertices. Otherwise, any convex body could be approximated arbitrarily well by a polytope in the class. Let $p\geq d+1$ be an integer and let $K$ have no more than $p$ vertices. In addition, suppose that the vertices of $K$ are not too peaked: Assume the existence of $\nu\in (0,1)$ and $\varepsilon_0>0$ such that
$\DS \textsf{Vol}_d\left(B_d(v,\varepsilon)\cap K\right)\geq \nu \textsf{Vol}_d\left(B_d(v,\varepsilon)\right)=\nu\varepsilon^d\kappa_d$, for all $v\in V$ and $\varepsilon\in(0,\varepsilon_0)$.
Finally, also suppose that $K\subseteq [0,1]^d$. Then, $\textsf{Vol}_d(K)\leq 1$ and for $\varepsilon\in(0,\varepsilon_0)$, \eqref{PolytoCase} becomes
\begin{align} \label{PolytoCase2}
	\PP_K\left[\textsf{d}_{\textsf{H}}(K,\hat K_n)\geq \varepsilon\right] & \leq p\exp\left(-\nu\kappa_d n\varepsilon^d\right).
\end{align}

Denote by $\mathcal P(p,\nu,\varepsilon_0)$ the class of all polytopes $K\subseteq [0,1]^d$ that have no more than $p$ vertices and such that $\displaystyle{\textsf{Vol}_d\left(B(v,\varepsilon)\cap K\right)\geq \nu\varepsilon^d\kappa_d}$, for all $v\in V$ and $\varepsilon\in(0,\varepsilon_0)$. This last condition is called a \textit{standardness condition}, see \cite{Cuevas2009}. We have the following theorem.

\begin{theorem}
	Let $d$ and $p\geq d+1$ be positive integers, $\nu\in (0,1)$ and $\varepsilon_0$ be positive real numbers. The following uniform deviation inequality holds:
	$$\sup_{K\in\mathcal P(p,\nu,\varepsilon_0)} \PP_K\left[n^{1/d}\textsf{d}_{\textsf{H}}(K,\hat K_n)\geq x\right] \leq pe^{-\nu\kappa_d x^d},$$
for all $x\in[0,\epsilon_0n^{1/d}]$.
\end{theorem}

Note that if $x\in (\varepsilon_0 n^{1/d}, \sqrt{d}n^{1/d})$, the right hand side can be replaced by $\DS pe^{-L\varepsilon_0^d n}$ and if $x>\sqrt{d}n^{1/d}$, it can be replaced by zero, since $\hat K_n\subseteq K\subseteq [0,1]^d$, yielding $\textsf{d}_{\textsf{H}}(\hat K_n,K)\leq\sqrt{d}$ almost surely. The next corollary follows directly. 

\begin{corollary} \label{CorPolSupp}
	Let $d$ and $p\geq d+1$ be positive integers, and $\nu\in (0,1)$ and $\varepsilon_0>0$. Then, for all real numbers $q\geq 1$,
	$$\sup_{K\in\mathcal P(p,\nu,\varepsilon_0)} \E_K\left[\textsf{d}_{\textsf{H}}(\hat K_n,K)^q\right]=O\left(n^{\frac{-q}{d}}\right).$$
\end{corollary}

It is easy to see from its proof that Theorem \ref{Thm:OptPol} still holds if the condition ``$B_d(a,\eta)\subseteq K$ for some $a\in\R^d$" is added under the supremum, for some small fixed $\eta<1$. Moreover, together with Proposition \ref{Claim:Gen} in Section \ref{Subsec:Spec}, Corollary \ref{corollary1} entails that 
	$$\sup_{K\in\mathcal K_d^{(1)}:\exists a\in\R^d,B_d(a,\eta)\subseteq K} \E_K\left[\textsf{d}_{\textsf{H}}(\tilde K_n,K)^q\right] = O\left(\left(\frac{\ln n}{n}\right)^{q/d}\right).$$
	The lower and upper bounds differ only by a logarithmic factor. Moreover, by adapting the proof of Theorem \ref{Thm:OptPol}, one can see that the rate in Corollary \ref{CorPolSupp} is tight and it is the same rate as in Theorem \ref{Thm:OptPol}. This yields that when $\mu$ is a uniform distribution, polytopal supports are the hardest convex bodies to estimate in a minimax sense, with respect to the Hausdorff metric. The contrary happens under the Nikodym metric: polytopal supports are the ones that can be estimated with the highest accuracy and convex supports with smooth boundary are the hardest to estimate, see \cite{Brunel2017}. In statistical terms, the estimator $\hat K_n$ does not adapt to polytopal supports under the Hausdorff metric, unlike under the Nikodym metric \cite{Brunel2016}. This is very surprising, since estimating $K$ under the Hausdorff metric is the same as estimating $h_K$ under the $L^\infty$ norm. Yet, when $K$ is a polytope, its support function is piecewise linear, so one would expect it to be easier to estimate.

\section{Uniform case: A rescaled version of the Hausdorff metric} \label{Sec:Resc}

In \cite{Brunel2016} and \cite{Brunel2017}, the risk of $\hat K_n$ is measured using the Nikodym distance rescaled by the volume of $K$, namely, $\displaystyle{\frac{\textsf{Vol}_d\left(\hat K_n\triangle K\right)}{\textsf{Vol}_d(K)}}$. Rescaling with the volume of $K$ allows to consider convex bodies $K$ that are not necessarily included in some given bounded set, e.g., $B_d(0,1)$. In this section, we introduce a scaled version of $\textsf{d}_{\textsf{H}}$ which allows to remove the assumption $K\subseteq B_d(0,1)$. In order to avoid unnecessary technicalities, we only deal with uniform distributions. 

Define the pseudo-distance $\textsf{d}_{\textsf{L}}$ between two sets $K$ and $K'$ as $\DS \textsf{d}_{\textsf{L}}(K,K')=\inf\{\varepsilon\geq 0 :\exists a\in R^d, a+(1+\varepsilon)^{-1}(K-a)\subseteq K'\subseteq a+(1+\varepsilon)(K-a)\}$.
It is easy to see that $\DS \ln\left(1+\textsf{d}_{\textsf{L}}\right)$ is a proper distance. If $\textsf{d}_{\textsf{L}}(K,K')$ is small, then $\DS \ln\left(1+\textsf{d}_{\textsf{L}}(K,K')\right)\approx \textsf{d}_{\textsf{L}}(K,K')$ and it is simpler, in our next results, to deal with $\DS \textsf{d}_{\textsf{L}}(K,K')$. One advantage of $\textsf{d}_{\textsf{L}}$ over $\textsf{d}_{\textsf{H}}$ is that it is invariant under the action of any invertible affine map $T$: $\textsf{d}_{\textsf{L}}(K,K')=\textsf{d}_{\textsf{L}}(T(K),T(K'))$, for all compact sets $K, K'$. Moreover, $\textsf{d}_{\textsf{L}}$ the well studied Banach-Mazur multiplicative distance. The Banach-Mazur distance between two convex bodies $K$ and $K'$ is defined as $\DS \textsf{d}_{\textsf{BM}}(K,K')=\inf\{\lambda\geq 1: \exists u,v\in\R^d, \exists T\in\textsf{GL}(d), K+u\subseteq T(K'+v)\subseteq \lambda(K+u)\}$, where $\textsf{GL}(d)$ stands for the set of all linear isomorphisms in $\R^d$ (see, e.g., \cite{Schneider1993}). It is a multiplicative pseudo-distance, in the sense that $\ln \textsf{d}_{\textsf{BM}}$ satisfies the triangle inequality. The following lemma holds.
\begin{lemma} \label{LemmaBM}
	For all $K,K'\in\K$, $\DS \textsf{d}_{\textsf{BM}}(K,K')\leq \left(1+\textsf{d}_{\textsf{L}}(K,K')\right)^2$.
\end{lemma}

\begin{proof}	
	Let $K,K'\in\K$. Let $\varepsilon>\textsf{d}_{\textsf{L}}(K,K')$. Then, there exists $a\in\R^d$ for which $\DS (1+\varepsilon)^{-1}(K-a)\subseteq K'-a\subseteq (1+\varepsilon)(K-a)$, yielding $\DS \DS K-a\subseteq (1+\varepsilon)(K'-a)\subseteq (1+\varepsilon)^2(K-a)$, i.e., $\DS (1+\varepsilon)^2\geq \textsf{d}_{\textsf{BM}}(K,K')$. Lemma \ref{LemmaBM} follows by letting $\varepsilon$ go to $\textsf{d}_{\textsf{L}}(K,K')$.	
\end{proof}

\noindent In particular, $\ln \textsf{d}_{\textsf{BM}}(K,K')\leq 2\textsf{d}_{\textsf{L}}(K,K')$, for all convex bodies $K, K'$. The following lemma, where we denote by $\textsf{conv}$ the convex hull, is straightforward. 

\begin{lemma} \label{Lemma:Aff}
	Let $T:\R^d\rightarrow\R^d$ be an invertible affine map. Let $K\in\K$, $X_1,\ldots,X_n$ be i.i.d. uniform random variables in $K$ and $Y_1,\ldots,Y_n$ be i.i.d. uniform random variables in $TK$. Then, $\textsf{d}_{\textsf{L}}(K,\textsf{conv}(X_1,\ldots,X_n))$ and $\textsf{d}_{\textsf{L}}(TK,\textsf{conv}(Y_1,\ldots,Y_n))$ have the same distribution.
\end{lemma}

As a consequence, working with $\textsf{d}_{\textsf{L}}$ instead of $\textsf{d}_{\textsf{H}}$ allows to remove the assumption $K\subseteq B_d(0,1)$.

\subsection{Convex supports with positive reach complement}

For $0<r\leq 1$, denote by $\tilde {\mathcal K}_{d,r}$ the class of all convex bodies $K\in\K$ satisfying the $\tilde r$-rolling condition, with $\displaystyle{\tilde r=r\left(\frac{\textsf{Vol}_d(K)}{\kappa_d}\right)^{1/d}}$. 
By introducing the number $\tilde r$ in this definition, we ensure that all convex bodies that are similar to $K$, i.e., obtained by rescaling, translating, rotating or reflecting $K$, will be in $\tilde {\mathcal K}_{d,r}$ as soon as $K$ is.

\begin{theorem}\label{theorem1}
	Let $d$ and $n$ be positive integers and $r\in (0,1]$. Set $\DS c=\frac{2d}{(d+1)\tau}$, $\displaystyle{\tau_3=\max\left(1,\frac{2^{(d+1)/2}d\kappa_d}{r^d\kappa_{d-1}}\right)}$, $\displaystyle{a_n=\left(\frac{\tau_3\ln n}{n}\right)^\frac{2}{d+1}}$ and $\displaystyle{b_n=n^{\frac{-2}{d+1}}}$. Then, the random polytope $\hat K_n$ satisfies
	$$\sup_{K\in\tilde{\mathcal K}_{d,r}} \PP_K\left[\textsf{d}_{\textsf{L}}(\hat K_n,K)\geq 4a_n + 4b_n x\right] \leq  6^d\exp\left(-c x^\frac{d+1}{2}\right),$$
for all $x\geq 0$ with $a_n+b_n x\leq 1/4$.
\end{theorem}

\begin{proof}

The proof of this theorem follows the same lines as that of Theorem \ref{theorem21}. Let $d\geq 1$, $r\in (0,1]$ and $K\in\tilde{\mathcal K}_{d,r}$. Let $\DS \tilde r=r\left(\frac{\textsf{Vol}_d(K)}{\kappa_d}\right)^{1/d}$. Let $n$ be a positive integer and $X, X_1, X_2, \ldots, X_n$ be $n+1$ i.i.d. random points uniformly distributed in $K$. Without loss of generality, one can assume that $0$ is an interior point of $K$ and $B_d(0,\tilde r)\subseteq K$. Indeed, since $K\in\tilde{\mathcal K}_{d,r}$, there exists a point $a\in K$ such that $B_d(a,\tilde r)\subseteq K$ and $K$ could be replaced with $K-a$ without changing the distribution of $\textsf{d}_{\textsf{L}}(\hat K_n,K)$. For all $u\in \Sp^{d-1}, h_K(u)\geq \tilde r$ and
\begin{align} \label{Distd_L}
	\textsf{d}_{\textsf{L}}(\hat K_n,K) & \leq \inf\{\varepsilon>0: (1+\varepsilon)^{-1}K\subseteq \hat K_n\subseteq (1+\varepsilon)K\} \nonumber \\
	& = \inf\{\varepsilon>0: (1+\varepsilon)^{-1}K\subseteq \hat K_n\} \nonumber \\
	& = \inf\{\varepsilon>0: \forall u\in \Sp^{d-1}, (1+\varepsilon)^{-1}h_K(u)\subseteq h_{\hat K_n}(u)\}.
\end{align}


Let $\varepsilon\in (0,1/4]$ and $u\in \Sp^{d-1}$. By the $\tilde r$-rolling condition, we have
\begin{align*}
	\textsf{Vol}_d\left(C_K(u,\varepsilon \tilde r/2)\right) & = \int_0^{\varepsilon \tilde r/2} \sqrt{x(2\tilde r-x)}^{d-1}\kappa_{d-1}dx \geq \kappa_{d-1}\tilde r^\frac{d-1}{2}\int_0^{\varepsilon \tilde r/2} x^\frac{d-1}{2}dx  \\
	& \geq \frac{2\kappa_{d-1}\tilde r^d}{2^{\frac{d+1}{2}}(d+1)}\varepsilon^\frac{d+1}{2} = \frac{2\kappa_{d-1} r^d \textsf{Vol}_d(K)}{2^{\frac{d+1}{2}}\kappa_d(d+1)}\varepsilon^\frac{d+1}{2}.
\end{align*}
Hence, if $\DS \eta=\frac{\kappa_{d-1} r^d}{2^{\frac{d-1}{2}}\kappa_d(d+1)}$,
\begin{align} \label{step113}
	\PP\left[\frac{h_{\hat K_n}(u)}{h_K(u)}\leq (1+\varepsilon)^{-1}\right] & \leq \PP\left[\frac{h_{\hat K_n}(u)}{h_K(u)}\leq 1-\frac{\varepsilon}{2}\right] \nonumber \\
	& = \PP\left[h_{\hat K_n}(u) \leq h_K(u)-\frac{\varepsilon h_K(u)}{2}\right]  \nonumber \\
	& \leq \PP\left[h_{\hat K_n}(u) \leq h_K(u)-\frac{\varepsilon \tilde r}{2}\right] \nonumber \\
	& = \left(1-\frac{\textsf{Vol}_d\left(C_K(u,\varepsilon \tilde r/2)\right)}{\textsf{Vol}_d(K)}\right)^n \nonumber \\
	& \leq \exp\left(- \eta  n\varepsilon^\frac{d+1}{2}\right).
\end{align}
By positive homogeneity of support functions (see \eqref{prop1}), \eqref{step113} holds for all $u\in\R^d\backslash\{0\}$. 

Let $\delta=\varepsilon/2$. Recall that $K^\circ\in\K$, since we have assumed that $0\in\textsf{int}(K)$, and that $0\in \textsf{int}(K^\circ)$ since $K$ is bounded. It is easy to adapt Lemma \ref{lemma1} and prove the existence of $\mathcal N_\delta\subseteq \partial(K^\circ)$, with cardinality at most $(3/\delta)^d$, such that every $u\in\partial(K^\circ)$ can be written as $\DS u=u_0-\sum_{j=1}^\infty \delta_j u_j$ with $0\leq \delta_j\leq \delta^j$ and $u_j\in\mathcal N_\delta$, for all $j\geq 1$. Then, using \eqref{step113} together with the union bound yields that with probability at least $\displaystyle{1-\left(\frac{6}{\varepsilon}\right)^d\exp\left(- \eta n\varepsilon^\frac{d+1}{2}\right)}$, 
\begin{equation} \label{event1}
	\frac{h_{\hat K_n}(u)}{h_K(u)} > (1+\varepsilon)^{-1}
\end{equation}
simultaneously for all $\DS u\in\mathcal N_\delta$. Assume that \eqref{event1} holds and let $u\in\partial (K^\circ)$. Write $\displaystyle{u=u_0-\sum_{j=1}^\infty \varepsilon_j u_j}$, as above. By \eqref{prop1}, \eqref{prop2} and 
\eqref{prop3}, 
\begin{equation} \label{step114}
	h_{\hat K_n}(u)\geq h_{\hat K_n}(u_0)-\sum_{j=1}^\infty \delta_j h_{\hat K_n}(u_j).
\end{equation}
In addition, for all $j\geq 1$, $(1+\varepsilon)^{-1}h_K(u_j)\leq h_{\hat K_n}(u_j) \leq h_K(u_j)$. Since $0\in \textsf{int}(K)$ and $u_j\in\partial (K^\circ)$, $h_K(u_j)=\|u_j\|_{K^\circ}=1$ and the following is true, for all $u\in\partial (K^\circ)$:
\begin{align} \label{step115}
	h_{\hat K_n}(u) & > (1+\varepsilon)^{-1} -\sum_{j=1}^\infty \delta^j \nonumber \\
	& \geq (1+\varepsilon)^{-1}-\frac{\delta}{1-\delta} \geq 1-2\varepsilon.
\end{align}
Since $0\in \textsf{int}(K^\circ)$, for all $u\in \Sp^{d-1}$, there exists $\lambda>0$ for which $\lambda u\in\partial (K^\circ)$, yielding 
\begin{equation*}
	\frac{h_{\hat K_n}(u)}{h_K(u)}=\frac{h_{\hat K_n}(\lambda u)}{h_K(\lambda u)} > 1-2\varepsilon\geq (1+4\varepsilon)^{-1},
\end{equation*}
since $\varepsilon\leq 1/4$.
Hence, with probability at least $\displaystyle{1-\left(\frac{6}{\varepsilon}\right)^d\exp\left(- \eta n\varepsilon^\frac{d+1}{2}\right)}$, $\DS \frac{h_{\hat K_n}(u)}{h_K(u)} > (1+4\varepsilon)^{-1}$
simultaneously for all $u\in \Sp^{d-1}$, yielding 
$\DS \textsf{d}_{\textsf{L}}(\hat K_n,K)< 4\varepsilon$.
Therefore, $\DS \PP\left[\textsf{d}_{\textsf{L}}(\hat K_n,K)\geq 4\varepsilon \right] \leq 6^d\exp\left(- \eta n\varepsilon^\frac{d+1}{2}-d\ln\varepsilon\right)$, which yields Theorem \ref{theorem1} by setting $\varepsilon=a_n+b_n x$.

\end{proof}

%
%
%

\subsection{General convex supports}

Now, we consider any $K\in\K$ and let $\mu$ be the uniform probability measure in $K$. By John's ellipsoid theorem (see \cite{Leichtweiss1959} for instance), there exists an invertible affine map $T:\R^d\rightarrow\R^d$ and $a\in\R^d$ such that $B_d(a,1/(2d))\subseteq TK\subseteq B_d(0,1/2)$. Necessarily, $a\in B_d(0,1/2)$, hence, up to a translation, one can assume that $B_d(0,1/(2d))\subseteq TK\subseteq B_d(0,1)$ without loss of generality. Therefore, by Lemma \ref{Lemma:Aff}, we can even assume that $B_d(0,1/(2d))\subseteq K\subseteq B_d(0,1)$, without loss of generality. Hence, using Proposition \ref{Claim:Gen}, one can adapt Theorem \ref{theorem21} with $\alpha=d$ in order to prove the following result.

\begin{theorem} \label{Thm:ScaledGen}
	Let $\DS a_n=\left(\frac{\ln n}{n}\right)^{1/d}$, $b_n=n^{-1/d}$ and $\DS L= \frac{(2d)^{d-1}\kappa_{d-1}}{3\kappa_d}$. Then, the random polytope $\hat K_n$ satisfies
	$$\sup_{K\in\K} \PP_K\left[\textsf{d}_{\textsf{L}}(\hat K_n,K)\geq 8da_n + 8db_n x\right] \leq  12^d\exp\left(-Lx^d\right),$$
for all $x\geq 0$ with $a_n+b_nx\leq 1/(16d)$.
\end{theorem}

\begin{proof}
Let $K\in\K$. As mentioned above, one can assume $B_d(0,1/(2d))\subseteq K\subseteq B_d(0,1)$, without loss of generality. By \eqref{Distd_L}, for all $\varepsilon\in (0,1/2]$, 

\begin{align*}
	\PP_K[\textsf{d}_{\textsf{L}}(\hat K_n,K)\geq \varepsilon] & \leq \PP_K\left[\exists u\in\Sp^{d-1}, h_{\hat K_n}\leq (1+\varepsilon)^{-1} h_K(u)\right] \\
	& \leq \PP_K\left[\exists u\in\Sp^{d-1}, h_{\hat K_n}\leq h_K(u)-\varepsilon/(4d)\right] \\
	& \leq \PP_K\left[\textsf{d}_{\textsf{H}}(K,\hat K_n)\geq \varepsilon/(4d)\right].
\end{align*}
Take $\varepsilon=8d(a_n+b_n x)$, where $a_n$ and $b_n$ are defined in the statement of Theorem \ref{Thm:ScaledGen}. Then, the conclusion follows directly from Theorem \ref{theorem21} and Proposition \ref{Claim:Gen} with $\eta=1/(2d)$, $\alpha=d$, $L= \frac{(2d)^{d-1}\kappa_{d-1}}{3\kappa_d}$ and by noting that $L\geq 1$ (so $\tau_1=1$), by \cite[Lemma 2.2]{ReisnerSchuttWerner2001}.

\end{proof}

%
%

\section{Random polytopes in high dimension} \label{Sec:CompComp}

In high dimension, the random polytope $\hat K_n$ is hard to compute. Indeed, it requires $O(n^{d/2})$ operations (see, e.g., \cite{Chazelle1993}), which becomes very costly when both $n$ and $d$ are large. The random polytope can be written using infinitely many constraints: $\DS \hat K_n=\{x\in\R^d:\langle u,x\rangle\leq \max_{1\leq i\leq n}\langle u,X_i\rangle, \forall u\in\Sp^{d-1}\}$. In this section, we show that if the parameter $\alpha$ in Assumption \ref{Assump} is known, then it is possible to significantly reduce the computational cost of $\hat K_n$ by taking only a finite subset of the linear constraints that define $\hat K_n$.

For an integer $M\geq 1$, let $U_1,\ldots,U_M$ be i.i.d. random unit vectors and let $$\hat K_{n,M}= \{x\in\R^d:\langle U_j,x\rangle\leq \max_{1\leq i\leq n}\langle U_j,X_i\rangle, \forall j=1,\ldots,M\}.$$
Then, the following holds.

\begin{theorem} \label{ThmCompComp}

Let $K\in\K$ and $X_1,\ldots,X_n$ be i.i.d. random points in $K$, with some probability measure $\mu$. Let Assumption \ref{Assump} hold. Set $M=Cn^{\frac{d-1}{\alpha}}(\ln n)^{\frac{\alpha-d+1}{\alpha}}$, where $C>2\alpha^{-1}d^2 8^{\frac{d-1}{2}}$. Then,
$$\textsf{d}_{\textsf{H}}(\hat K_{n,M},K)=O_{\PP}\left(\left(\frac{\ln n}{n}\right)^{\frac{1}{\alpha}}\right).$$

\end{theorem}

\begin{proof}

For $j=1,\ldots,M$, let $\hat h_j=\max_{1\leq i\leq n}\langle U_j,X_i\rangle=h_{\hat K_n}(U_j)$ and $h_j=h_K(U_j)$. 
Since $K\in\K$, $K$ has nonempty interior and it is bounded. Thus, let $R>r>0$ and $a\in\R^d$ such that $B(a,r)\subseteq G\subseteq B(a,R)$. 
Let $\eta\in (0,r)$ and $\varepsilon\in (0,1)$. Let $\mathcal A$ be the event: ``$|\hat h_j-h_j|\leq \eta$, for all $j=1,\ldots,M$" and $\mathcal B$ be the event ``The collection $\{U_1,\ldots,U_M\}$ is an $\varepsilon$-net of $\Sp^{d-1}$".

Using \eqref{step1} together with the union bound,

\begin{equation} \label{ProbaAcomp}
	\PP[\mathcal A^{\complement}]\leq Me^{-Ln\eta^\alpha}.
\end{equation}

By \cite[Lemma10]{BrunelTDLS2017},
 
\begin{equation} \label{ProbaNets}
	\PP[\mathcal B^{\complement}]\leq 6^d\exp\left(-\frac{M\varepsilon^{d-1}}{2d8^{\frac{d-1}{2}}}+d\ln\left(\frac{1}{\varepsilon}\right)\right).
\end{equation}

Assume that both $\mathcal A$ and $\mathcal B$ are satisfied. Then, by \cite[Lemma 7]{BrunelTDLS2017}, $\DS \textsf{d}_{\textsf{H}}(\hat K_n,K)\leq \frac{\eta R}{r}\frac{1+\eta/r}{1-\eta/r}+\frac{2R\varepsilon}{1-\varepsilon}$. As a consequence of \eqref{ProbaAcomp} and \eqref{ProbaNets}, this yields 

\begin{equation}
	\PP[\textsf{d}_{\textsf{H}}(\hat K_{n,M},K)> \frac{\eta R}{r}\frac{1+\eta/r}{1-\eta/r}+\frac{2R\varepsilon}{1-\varepsilon}] \leq Me^{-Ln\eta^\alpha}+ 6^d\exp\left(-\frac{M\varepsilon^{d-1}}{2d8^{\frac{d-1}{2}}}+d\ln\left(\frac{1}{\varepsilon}\right)\right).
\end{equation}

Setting $\DS \varepsilon=\left(\frac{\ln n}{n}\right)^{1/\alpha}$ and $\DS \eta=\left(\frac{x\ln n}{n}\right)^{1/\alpha}$, for some positive number $x$, ensuring that $\eta<r$, completes the proof of Theorem \ref{ThmCompComp}.

\end{proof}

As a corollary, the random polytope $\hat K_n$ can be approximated by $\hat K_{n,M}$:

\begin{corollary}\label{CorApproxRP}

Let the assumptions of Theorem \ref{ThmCompComp} hold. Then,
$$\textsf{d}_{\textsf{H}}(\hat K_{n,M},K)=O_{\PP}\left(\left(\frac{\ln n}{n}\right)^{\frac{1}{\alpha}}\right).$$

\end{corollary}

Computing $\hat K_{n,M}$ requires sampling $M$ uniform random unit vectors (which can be done by sampling standard Gaussian vectors and rescaling them) and checking $nM$ linear inequalities. Hence, the overall computation cost of computing $\hat K_{n,M}$ is $O(nMd)$. For instance, when $\mu$ is the uniform distribution in a general convex body $K$, we saw previously that $\alpha$ can be chosen as $\alpha=d$. Hence, Corollary \ref{CorApproxRP} shows that $\hat K_n$ can approximated in $\DS O\left(d^28^{\frac{d-1}{2}}n^{2-1/d}(\ln n)^{1/d}\right)$ operations, without affecting the rate of its statistical accuracy. If $\mu$ is the uniform distribution in a convex body whose complement has positive reach, then $\hat K_n$ can be approximated in $\DS O\left(d^28^{\frac{d-1}{2}}n^{\frac{3d-1}{d+1}}(\ln n)^{-\frac{d-3}{d+1}}\right)$ operations. In both cases, the computation cost is significantly better than that of $\hat K_n$, which is $O(n^{d/2})$. Note that the computation cost of $\hat K_{n,M}$ is still exponential in the dimension, through the multiplicative constants.

\bibliographystyle{plain}
\bibliography{Biblio}

\end{document}